\documentclass[11pt]{amsart}
\usepackage{amssymb, latexsym}
\theoremstyle{plain}
\newtheorem{theorem}{Theorem}

\newtheorem{proposition}{Proposition}

\newtheorem*{1'}{Theorem 1-Bessel}
\newtheorem*{P2'}{Proposition 2-Bessel}
\newtheorem*{P3'}{Proposition 3-Bessel}
\newtheorem*{P4'}{Proposition 4-Bessel}
\newtheorem*{C1'}{Corollary 1-Bessel}

\newtheorem*{2'}{Theorem 2-Bessel}
\newtheorem*{3'}{Theorem 3-Bessel}

\theoremstyle{remark}

\newtheorem*{Remark 1}{Remark 1}
\newtheorem*{Remark 2}{Remark 2}
\newtheorem*{Remark 3}{Remark 3}
\newtheorem*{Remark 4}{Remark 4}

\numberwithin{equation}{section}
\renewcommand{\baselinestretch}{1.4}

\begin{document}

\title[Large deviations for longest alt/inc subsequence] {Large deviations for the longest alternating and the longest increasing subsequence in a random  permutation avoiding a pattern of length three}

\author{Ross G. Pinsky}


\address{Department of Mathematics\\
Technion---Israel Institute of Technology\\
Haifa, 32000\\ Israel}
\email{ pinsky@technion.ac.il}
\urladdr{https://pinsky.net.technion.ac.il/}

\date{}

\begin{abstract}
We calculate the large deviations for the length of the longest alternating subsequence and for the length of the longest increasing subsequence
in a uniformly random permutation that avoids a pattern of length three.  We treat all six patterns in the case of alternating subsequences. 
In the case of increasing subsequences, we treat two of the three patterns for which a classical large deviations result is possible. 
The same rate function appears in all six cases for alternating subsequences. This rate function is in fact the rate function 
for the large deviations of the sum of IID symmetric Bernoulli random variables. The same rate function appears in the two cases we treat for increasing
subsequences. This rate function is twice the rate function for alternating subsequences.

\end{abstract}

\maketitle
\section{Introduction and Statement of Results}\label{intro}
\renewcommand{\baselinestretch}{1.3}

The problem of analyzing  the distribution of the length, $L_n$, of the longest
increasing subsequence  in a uniformly random permutation
from $S_n$, the set of permutations of $[n]:=\{1,\cdots, n\}$,
has a long and
distinguished history; see \cite{AD} and references therein, and see \cite{R}. In
particular, the work of  Logan and Shepp
\cite{LS} together with that   of Vershik and Kerov \cite{VK}
show that $EL_n\sim 2n^\frac12$
as $n\to\infty$.
This was followed over twenty years later  by the  profound work of   Baik, Deift and Johansson \cite{BDJ}, who proved
that
$$
\lim_{n\to\infty}P(\frac{L_n-2n^\frac12}{n^\frac16}\le x)=
F(x),
$$
where  $F$ is the Tracy-Widom distribution.
A large deviations result for the lower tail probabilities $P(\frac{L_n}{n^\frac12}\le x)$, for $x<2$, was given in
\cite{DZ99}, while  for the upper tail probabilities  $P(\frac{L_n}{n^\frac12}\ge x)$, for $x>2$, one was given in \cite{S}.  See also references therein.

Now consider the length of the longest alternating subsequence in a uniformly random permutation from $S_n$. An alternating subsequence of length $k$ in a permutation $\sigma=\sigma_1\cdots\sigma_n\in S_n$ is a subsequence of the form
$\sigma_{i_1}>\sigma_{i_2}<\sigma_{i_3}>\sigma_{i_4}\cdots\sigma_{i_k}$
or $\sigma_{i_1}<\sigma_{i_2}>\sigma_{i_3}<\sigma_{i_4}\cdots\sigma_{i_k}$,
 where $1\le i_1<\cdots<i_k\le n$.
Call the first type an initially descending alternating subsequence and call the second type an initially ascending alternating subsequence.
Of course, for asymptotic results concerning the longest alternating subsequence, it doesn't matter which type one considers since the two differ from one another by at most one.
Stanley derived the exact expected value and variance for initially descending alternating subsequences \cite{St}.
In particular, letting $A_n$ denote the length of the longest alternating subsequence of either type, he showed
that $E_nA_n\sim\frac23 n$ and that  the variance is asymptotic to $\frac8{45}n$. It has been proven that  a central limit theorem holds with a Gaussian limiting distribution
\cite{W,St}. We are unaware of large deviations results for this permutation statistic.
\medskip

In this paper we derive the large deviations for increasing subsequences and alternating subsequences in uniformly  random permutations that avoid a particular pattern in $S_3$.
We recall the definition of a pattern avoiding permutation.
If $\sigma=\sigma_1\cdots\sigma_n\in S_n$ and $\eta=\eta_1\cdots\eta_m\in S_m$, where $2\le m\le n$,
then we say that $\sigma$ contains $\eta$ as a pattern if there exists a subsequence $1\le i_1<i_2<\cdots<i_m\le n$ such
that for all $1\le j,k\le m$, the inequality $\sigma_{i_j}<\sigma_{i_k}$ holds if and only if the inequality $\eta_j<\eta_k$ holds.
If $\sigma$ does not contain $\eta$, then we say that $\sigma$ \it avoids\rm\ $\eta$.
We denote by $S_n^{\text{av}(\eta)}$ the set of permutations in $S_n$ that avoid $\eta$.
For any $\eta\in S_m$, we denote the uniform probability measure on $S_n^{\text{av}(\eta)}$ by $P_n^{\text{av}(\eta)}$ and denote expectations by
$E_n^{\text{av}(\eta)}$.

It is well-known \cite{B,SS} that $|S_n^{\text{av}(\eta)}|=C_n$, for every $\eta\in S_3$, where $C_n=\frac1{n+1}\binom{2n}{n}$
is the $n$th Catalan number. One has
\begin{equation}\label{Catasymp}
C_n\sim\frac{4^n}{\sqrt\pi n^\frac32}.
\end{equation}

We begin with alternating subsequences.
In \cite{FMW} it was proven that   $E_n^{\text{av}(\eta)}A_n\sim \frac n2$,
that the variance  of $A_n$ under
$P_n^{\text{av}}(\eta)$
is asymptotic to $\frac14n$
and  that $\frac{A_n-\frac n2}{\frac12\sqrt n}$ converges in distribution to the standard Gaussian distribution, for all choices of $\eta\in S_3$.
We will prove the following theorem.
\begin{theorem}\label{altthm}
Let $\eta\in S_3$. The longest alternating subsequence $A_n$ satisfies
\begin{equation}\label{ldalt}
\begin{aligned}
&\lim_{n\to\infty}\frac1n\log P_n^{\text{av}(\eta)}(A_n\ge nx)=\lim_{n\to\infty}\frac1n\log P_n^{\text{av}(\eta)}(A_n> nx)=
-I^{\text{alt}}(x),\ x\in[\frac12,1);\\
&\lim_{n\to\infty}\frac1n\log P_n^{\text{av}(\eta)}(A_n\le nx)=\lim_{n\to\infty}\frac1n\log P_n^{\text{av}(\eta)}(A_n< nx)=-I^{\text{alt}}(x),\ x\in(0,\frac12],
\end{aligned}
\end{equation}
where
\begin{equation}\label{I-alt}
I^{\text{alt}}(x)=x\log x+(1-x)\log(1-x)+\log 2,\ x\in(0,1).
\end{equation}
\end{theorem}
\noindent \bf Remark 1.\rm\ The function $I^{\text{alt}}$ above is in fact the relative entropy $H(\mu_x;\mu_\frac12)$ of $\mu_x$ with respect to $\mu_\frac12$, where $\mu_p$ 
denotes the 
distribution of the Bernoulli random variable $X_p$ satisfying $P(X_p=1)=1-P(X_p=0)=p$.
 From Cram\`ers theorem, if $\{X_n\}_{n=1}^\infty$ are IID Bernoulli random variables with parameter $\frac12$, and $S_n=\sum_{j=1}^nX_j$, then
\eqref{ldalt} also holds with
$P_n^{\text{av}(\eta)}(A_n\cdots)$ replaced by $P(S_n\cdots)$.
in all four places.
\it It would be very interesting to understand why this connection arises.\rm\

\noindent \bf Remark 2.\rm\  Since $|S_n^{\text{av}(\eta)}|=C_n$, it follows
from \eqref{Catasymp} and \eqref{ldalt} that
$$
\begin{aligned}
&\lim_{\epsilon\to0^+}\lim_{n\to\infty}\frac1n\log |\{\sigma\in S_n^{\text{av}(\eta)}: A_n(\sigma)\ge1-\epsilon\}|=\\
&\lim_{\epsilon\to0^+}\lim_{n\to\infty}\frac1n\log |\{\sigma\in S_n^{\text{av}(\eta)}: A_n(\sigma)\le\epsilon\}|=
\log 2.
\end{aligned}
$$
\medskip

We now turn to  increasing subsequences.
In \cite{DHW}, the asymptotic behavior of the expectation $E_n^{\text{av}(\eta)}L_n$  and the variance $v_n(\eta)$ of the longest increasing subsequence $L_n$ under
$P_n^{\text{av}(\eta)}$ were obtained for all six permutations $\eta\in S_3$. Of course, the case $\eta=123$ is trivial. The
 expectation $E_n^{\text{av}(\eta)}L_n$ is on the order $n$ only for $\eta\in\{231,312,321\}$.
The limiting distribution
of $\frac{L_n-E_n^{\text{av}(\eta)}L_n}{v_n(\eta)}$ was calculated as well, the limit being
Gaussian only for  $\eta\in\{231,312\}$.
The paper culled a lot of other results in the literature in order to proceed.
From the results, it is clear that
a classical large deviations result is only possible for $\eta\in\{231,312,321\}$.
We consider here
 $\eta\in\{231,312\}$. In both of these cases, it was shown in  \cite{DHW} that $E_n^{\text{av}(\eta)}L_n=\frac{n+1}2$.
We will prove the following theorem.
\begin{theorem}\label{incthm}
Let $\eta\in\{231,312\}$. The longest increasing subsequence $L_n$ satisfies
\begin{equation}\label{ldinc}
\begin{aligned}.
&\lim_{n\to\infty}\frac1n\log P_n^{\text{av}(\eta)}(L_n\ge nx)=-I^{\text{inc}}(x),\ x\in[\frac12,1];\\
&\lim_{n\to\infty}\frac1n\log P_n^{\text{av}(\eta)}(L_n\le nx)=-I^{\text{inc}}(x),\ x\in(0,\frac12],\\
\end{aligned}
\end{equation}
where
\begin{equation}\label{I-inc}
\begin{aligned}
&I^{\text{inc}}(x)=2I^{\text{alt}}(x)=2\big(x\log x+(1-x)\log(1-x)+\log2\big), \ x\in(0,1);\\
&I^{\text{inc}}(1)=\log4.
\end{aligned}
\end{equation}
\end{theorem}
\bf \noindent Remark 1.\rm\ 
The function $I^{\text{inc}}$ above is in fact the relative entropy $H(\nu_x;\nu_\frac12)$ of $\nu_x$ with respect to $\nu_\frac12$, where $\nu_p$
denotes the
distribution of one-half the sum of two independent Bernoulli random variables with parameter $p$. 
 From Cram\`ers theorem, if $\{Y_n\}_{n=1}^\infty$ are IID random variables with distribution $\frac12(X^{(1)}_\frac12+X^{(2)}_\frac12)$, where
 $X^{(1)}_\frac12$ and $X^{(2)}_\frac12$ are IID Bernoulli random variables with parameter $\frac12$, and $S_n=\sum_{j=1}^n Y_j$, then
\eqref{ldinc} also holds with
$P_n^{\text{av}(\eta)}(L_n\cdots)$ replaced by $P(S_n\cdots)$.
in all four places.

\bf\noindent Remark 2.\rm\
The identity permutation in $S_n$ is the only permutation $\sigma\in S_n$ for which $L_n(\sigma)=n$. From this fact and \eqref{Catasymp} alone, it follows that
  $\lim_{n\to\infty}\frac1n\log P_n^{\text{av}(\eta)}(L_n\ge n)=-\log 4$.

\bf\noindent Remark 3.\rm\ The proof of Theorem \ref{incthm} uses generating functions. The same type of proof could be used to obtain the expectation and the variance
of $L_n$, which is considerably simpler than the proofs of these results  in \cite{DHW}.

The proof of Theorem \ref{altthm} is given in section \ref{pfaltthm} and the proof of Theorem \ref{incthm} in given in section \ref{pfincthm}.

\section{Proof of Theorem \ref{altthm}}\label{pfaltthm}
Recall that the \it reverse\rm\  of a permutation $\sigma=\sigma_1\cdots\sigma_n$ is the permutation $\sigma^{\text{rev}}:=\sigma_n\cdots\sigma_1$,
and the \it complement\rm\ of $\sigma$ is the permutation
$\sigma^{\text{com}}$ satisfying
$\sigma^{\text{com}}_i=n+1-\sigma_i,\ i=1,\cdots, n$.
Let $\sigma^{\text{rev-com}}$ denote the permutation obtained by applying reversal and then complementation to $\sigma$ (or equivalently, vice versa).
Although there are six permutations $\eta$ in $S_3$, to prove the theorem, it suffices to consider just two of them---one from $\{231,213,312, 132\}$ and one from $\{123,321\}$.
Indeed, for all $n\ge3$, the operation reversal is a bijection from $S_n^{\text{av}(231)}$ to $S_n^{\text{av}(132)}$ and from $S_n^{\text{av}(123)}$ to $S_n^{\text{av}(321)}$,
the operation complementation is a bijection from $S_n^{\text{av}(231)}$ to $S_n^{\text{av}(213)}$ and the operation reversal-complementation is a bijection from
$S_n^{\text{av}(231)}$
to $S_n^{\text{av}(312)}$. Furthermore,
 $A_n(\sigma)=A_n(\sigma^{\text{com}})=A_n(\sigma^{\text{rev-com}})=A_n(\sigma^{\text{rev-com}})$, for $\sigma\in S_n$. Proposition 2.2--iii in  \cite{FMW} shows that
the distributions of $A_n$ under $P_n^{\text{av}(231)}$  and under $P_n^{\text{av}(321)}$ coincide.
Thus, to prove the theorem, we need only consider the case $\eta=231$.

For $n\in\mathbb{N}$ and $\sigma\in S_n$, let $A_n^{+,-}(\sigma)$ denote the longest alternating subsequence in $\sigma$ that begins with an ascent and ends with a descent. An alternating subsequence that begins with an ascent and ends with of descent is of the form
$\sigma_{i_1}<\sigma_{i_2}>\sigma_{i_3}<\cdots>\sigma_{i_{2k+1}}$, for $1\le i_1<i_2<\cdots<i_{2k+1}\le n$, with $k\in\mathbb{N}$. If there is no such alternating subsequence
in $\sigma$, then define $A_n^{+,-}(\sigma)=1$.
Note that $A_n^{+,-}(\sigma)$ takes on positive, odd integral values.
It suffices to prove the theorem with $A_n^{+,-}(\sigma)$ in place of $A_n$ since $A_n(\sigma)-A_n^{+,-}(\sigma)\in\{0,1,2\}$, for all $\sigma\in S_n$.
For convenience in the proof, we define $A^{+,-}_0\equiv 0$.

Let
$$
M_n(\lambda)=E_n^{\text{av}(231)}e^{\lambda A^{+,-}_n},\ \lambda\in\mathbb{R}, n\ge 0,
$$
denote the moment generating function of $A^{+,-}_n$.
The main part of the proof of the theorem is the proof of the following proposition.
\begin{proposition}\label{generatingprop}
\begin{equation}\label{laplacetranslim}
\lim_{n\to\infty}\frac1n\log M_n(\lambda)=\log(e^\lambda+1)-\log 2.
\end{equation}
\end{proposition}
Let $I(x)$ denote the Legendre-Fenchel transform of the function appearing on the right hand side of \eqref{laplacetranslim}; that is,
\begin{equation}\label{I}
I(x)=\sup_{\lambda\in\mathbb{R}}\Big(\lambda x-\log(e^\lambda+1)+\log 2\Big),\ x\in\mathbb{R}.
\end{equation}
We have the following proposition.
\begin{proposition}\label{supL-F}
The function $I$,  defined in \eqref{I} and restricted to $x\in(0,1)$, is equal to $I^{\text{alt}}$ defined in \eqref{I-alt}.
\end{proposition}
Proposition \ref{supL-F} is well-known, as it corresponds to the case of IID symmetric Bernoulli random variables---see Remark 1 after Theorem \ref{altthm}; thus, we will skip the proof, which is a calculus exercise.
Using Propositions \ref{generatingprop} and \ref{supL-F}, the proof of Theorem \ref{altthm}
follows from a fundamental theorem in large deviations theory.

\noindent \it Proof of Theorem \ref{altthm}.\rm\ In light of Propositions \ref{generatingprop} and \ref{supL-F}, \eqref{ldalt} follows from the G\"artner-Ellis
theorem \cite{DZ}.
\hfill $\square$

It remains to prove  Proposition \ref{generatingprop}.

\noindent \it Proof of Proposition \ref{generatingprop}.\rm\
Every permutation $\sigma\in S_k^{\text{av}(231)}$  has the property that if $\sigma_j=n$, then the numbers $\{1,\cdots, j-1\}$ appear in the first $j-1$ positions in $\sigma$ (and then of course, the numbers $\{j,\cdots, n-1\}$ appear in the last $n-j$ positions in $\sigma$.)
From this fact, along with the fact that $|S_n^{\text{av}(\eta)}|=C_n$, it follows that
\begin{equation}\label{nprob}
P_n^{\text{av}(231)}(\sigma_j=n)=\frac{C_{j-1}C_{n-j}}{C_n},\ \text{for}\ j\in[n],
\end{equation}
where $C_0=1$.
It also follows that under the conditional measure $P_n^{\text{av}(231)}|\{\sigma_j=n\}$, the permutation $\sigma_1\cdots\sigma_{j-1}\in S_{j-1}$ has the distribution
$P_{j-1}^{\text{av}(231)}$,  the permutation $\sigma'_{j+1}\cdots\sigma'_n$ has the distribution $P_{n-j}^{\text{av}(231)}$, where $\sigma'_k=\sigma_k-j+1$, for $k=j+1,\cdots, n$, and these two permutations are independent.
From this last fact and the definition of $A_n^{+,-}$, it follows that
\begin{equation}\label{conddist}
\begin{aligned}
&A^{+,-}_n|\{\sigma_j=n\}\stackrel{\text{dist}}{=}A^{+,1}_{j-1}+1+A^{+,-}_{n-j},\ n\ge3,\ j=2,\cdots,n-1;\\
&A^{+,-}_n|\{\sigma_1=n\}\stackrel{\text{dist}}{=}A^{+,-}_n|\{\sigma_n=n\}\stackrel{\text{dist}}{=}A^{+,-}_{n-1},\ n\ge2,
\end{aligned}
\end{equation}
where on the right hand side of \eqref{conddist}, for any $k$, the random variable $A^{+,-}_k$ is considered on $S_k$ under the measure $P_k^{\text{av}(231)}$, and
where on the right hand side of  the first line in \eqref{conddist},  $A^{+,1}_{j-1}$ and $A^{+,-}_{n-j}$ are independent.

From \eqref{nprob} and \eqref{conddist}, we have
\begin{equation}\label{key}
\begin{aligned}
&M_n(\lambda)=\sum_{j=1}^n P_n^{\text{av}(231)}(\sigma_j=n)E_n^{\text{av}(231)}(e^{\lambda A^{+,-}_n}|\sigma_j=n)=\\
&\sum_{j=2}^{n-1}\frac{C_{j-1}C_{n-j}}{C_n}e^\lambda E^{\text{av}(231)}_{j-1}e^{\lambda A^{+,-}_{j-1}}E^{\text{av}(231)}_{n-j}e^{\lambda A^{+,-}_{n-j}}+2\frac{C_{n-1}}{C_n}
E^{\text{av}(231)}_{n-1}e^{\lambda A^{+,-}_{n-1}}=\\
&e^\lambda\sum_{j=2}^{n-1}\frac{C_{j-1}C_{n-j}}{C_n}M_{j-1}(\lambda)M_{n-j}(\lambda)+2\frac{C_{n-1}}{C_n}M_{n-1}(\lambda),\ n\ge 3.
\end{aligned}
\end{equation}
Multiplying the leftmost expression and the rightmost expression in \eqref{key} by
$C_nt^n$, we have
\begin{equation}\label{fromkey}
\begin{aligned}
&C_nM_n(\lambda)t^n=e^\lambda t\Big(\sum_{j=2}^{n-1}C_{j-1}M_{j-1}(\lambda)C_{n-j}M_{n-j}(\lambda)\Big)t^{n-1}+\\
&2tC_{n-1}M_{n-1}(\lambda)t^{n-1}, \ n\ge3.
\end{aligned}
\end{equation}
Summing over $n$ gives
\begin{equation}\label{fromkeywitht}
\begin{aligned}
&\sum_{n=3}^\infty C_nM_n(\lambda)t^n=e^\lambda t\sum_{n=3}^\infty\Big(\sum_{j=2}^{n-1}C_{j-1}M_{j-1}(\lambda)C_{n-j}M_{n-j}(\lambda)\Big)t^{n-1}+\\
&2t\sum_{n=3}^\infty C_{n-1}M_{n-1}(\lambda)t^{n-1}.
\end{aligned}
\end{equation}

Define
\begin{equation}\label{Gseries}
G_\lambda(t)=\sum_{n=0}^\infty C_nM_n(\lambda)t^n.
\end{equation}
For use in some of the  calculations  below, note that $C_0=C_1=1$, $C_2=2$, $M_0(\lambda)=1$ and $M_1(\lambda)=M_2(\lambda)=e^\lambda$.
We have
\begin{equation}\label{G^2}
G_\lambda^2(t)=\sum_{n=0}^\infty\Big(\sum_{k=0}^n C_kM_k(\lambda)C_{n-k}M_{n-k}(\lambda)\Big)t^n.
\end{equation}
The double sum on the right hand side of \eqref{fromkeywitht} can be written as
\begin{equation}\label{doublesum}
\begin{aligned}
&\sum_{n=3}^\infty\Big(\sum_{j=2}^{n-1}C_{j-1}M_{j-1}(\lambda)C_{n-j}M_{n-j}(\lambda)\Big)t^{n-1}=\\
&\sum_{n=3}^\infty\Big(\sum_{k=1}^{n-2}
C_kM_k(\lambda)C_{n-1-k}M_{n-1-k}(\lambda)\Big)t^{n-1}=\\
&\sum_{n=3}^\infty\Big(\sum_{k=0}^{n-1}
C_kM_k(\lambda)C_{n-1-k}M_{n-1-k}(\lambda)\Big)t^{n-1}-2\sum_{n=3}^\infty C_{n-1}M_{n-1}(\lambda)t^{n-1}=\\
&\sum_{m=2}^\infty\Big(\sum_{k=0}^m
C_kM_k(\lambda)C_{m-k}M_{m-k}(\lambda)\Big)t^m-2\sum_{m=2}^\infty C_mM_m(\lambda)t^m=\\
&(G^2_\lambda(t)-2e^\lambda t-1)-2(G_\lambda(t)-e^\lambda t-1)=G^2_\lambda(t)-2G_\lambda(t)+1,
\end{aligned}
\end{equation}
where \eqref{G^2} has been used for the penultimate  equality.
From \eqref{fromkeywitht}, \eqref{doublesum} and the definition of $G_\lambda$, we obtain
$$
\big(G_\lambda(t)-2e^\lambda t^2-e^\lambda t-1\big)
=e^\lambda t\big(G^2_\lambda(t)-2G_\lambda(t)+1\big)+2t\big(G_\lambda(t)-e^\lambda t-1\big),
$$
or equivalently
\begin{equation*}\label{Gequ}
e^\lambda tG^2_\lambda(t)+\big(2(1-e^\lambda)t-1\big)G_\lambda(t)+2(e^\lambda-1)t+1=0.
\end{equation*}
Since $G_\lambda(0)=1$, the quadratic formula yields
$$
G_\lambda(t)=\frac{1-2(1-e^\lambda)t-\sqrt{\big(2(1-e^\lambda)t-1\big)^2-4e^\lambda t\big(2(e^\lambda-1)t+1\big)}}{2e^\lambda t},
$$
which we rewrite as
\begin{equation}\label{G}
G_\lambda(t)=\frac{1-2(1-e^\lambda)t-\sqrt{4(1-e^{2\lambda})t^2-4t+1}}{2e^\lambda t}.
\end{equation}

From \eqref{G}, it follows that the radius of convergence $R_\lambda$ of the power series representing $G_\lambda$ is the smaller of the absolute values of the two roots
of the quadratic polynomial $4(1-e^{2\lambda})t^2-4t+1$.
The roots of this polynomial are $\frac{1\pm e^\lambda}{2(1-e^{2\lambda})}$.
Thus, we have $R_\lambda=\frac{|e^\lambda-1|}{2|e^{2\lambda}-1|}=\frac1{2(e^\lambda+1)}$.
Consequently,
from \eqref{Gseries}, we obtain
\begin{equation}\label{limsup}
\limsup_{n\to\infty}(C_nM_n(\lambda))^\frac1n=2(e^\lambda+1).
\end{equation}
From \eqref{Catasymp} and \eqref{limsup}, we obtain
\begin{equation}\label{limsupMn}
\limsup_{n\to\infty}\frac1n\log M_n(\lambda)=\log(e^\lambda+1)-\log 2.
\end{equation}
In light of \eqref{limsupMn}, to complete the proof of the proposition it suffices to show that the limit
$\lim_{n\to\infty}\frac1n\log M_n(\lambda)$ exists.

From  \eqref{Gseries}, \eqref{G} and \eqref{Catasymp}, in order to show that
$\lim_{n\to\infty}\frac1n\log M_n(\lambda)$ exists,
 it suffices to show that
$\lim_{n\to\infty}a_n^\frac1n$ exists, where
$a_n$ is the coefficient of $t^n$ in the power series expansion
\begin{equation}\label{a}
\sqrt{4(1-e^{2\lambda})t^2-4t+1}=\sum_{n=0}^\infty a_nt^n.
\end{equation}
Intuitively, it seems ``obvious'' that this limit exists, but unfortunately, we don't have a real quick  proof.
Recall that the two roots of $4(1-e^{2\lambda})t^2-4t+1$ are $r_1=\frac{1-e^\lambda}{2(1-e^{2\lambda})}>0$ and $r_2=\frac{1+e^\lambda}{2(1-e^{2\lambda})}$.
(We suppress the dependence on $\lambda$.)
We have $r_1<|r_2|$.
We write
\begin{equation}\label{polyroots}
4(1-e^{2\lambda})t^2-4t+1=(1-\frac t{r_1})(1-\frac t{r_2}).
\end{equation}
The Taylor series of $\sqrt{1-x}$
around $x=0$ is given by
\begin{equation}\label{Taylor}
\begin{aligned}
&\sqrt{1-x}=\sum_{n=0}^\infty (-1)^n\binom{\frac12}{n}x^n=1+\sum_{n=1}^\infty(-1)^n\frac{\frac12(\frac12-1)\cdots(\frac12-(n-1))}{n!}x^n=\\
&1-\sum_{n=1}^\infty\frac{(2n-2)!}{n!(n-1)!2^{2n-1}}x^n.
\end{aligned}
\end{equation}
Thus, from \eqref{polyroots} and \eqref{Taylor}, we have
\begin{equation}\label{explicitexpansion}
\begin{aligned}
&\sqrt{4(1-e^{2\lambda})t^2-4t+1}=\sqrt{1-\frac t{r_1}}\thinspace\sqrt{1-\frac t{r_2}}=\\
&1-\sum_{n=1}^\infty\frac{(2n-2)!}{n!(n-1)!2^{2n-1}r_1^n}t^n-
\sum_{n=1}^\infty\frac{(2n-2)!}{n!(n-1)2^{2n-1}!r_2^n}t^n\\
&+\sum_{n=2}^\infty\Big(\sum_{j=1}^n\frac{(2j-2)!}{j!(j-1)!2^{2j-1}\thinspace r_1^j}
\frac{(2(n-j)-2)!}{(n-j)!(n-j-1)!2^{2(n-j)-1}\thinspace r_2^{n-j}}\Big)t^n
\end{aligned}
\end{equation}
Using Stirling's formula, one finds that the expression $\frac{(2m-2)!}{m!(m-1)!2^{2m-1}}$, for $m\in\mathbb{N}$,  decays to zero on the order $m^{-\frac32}$. In particular then, this expression is bounded and has  sub-exponential
decay.
Using this with \eqref{a} and \eqref{explicitexpansion} and the fact that $r_1<|r_2|$,
it follows that $\lim_{n\to\infty}a_n^\frac1n=\frac1{r_1}$.
This completes the proof of Proposition \ref{generatingprop}\hfill$\square$

\section{Proof of Theorem \ref{incthm}}\label{pfincthm}
Consider $\eta\in\{231,312\}$. For convenience, let $L_0=0$.
Let
$$
M_n(\lambda)=E_n^{\text{av}(\eta)}e^{\lambda L_n},\ \lambda\in\mathbb{R},\ n\ge0,
$$
denote the moment generating function of $L_n$.
The main part of the proof of the theorem is the proof of the following proposition.
\begin{proposition}\label{generatingpropagain}
Let $\eta\in\{231,312\}$. Then
\begin{equation}\label{laplacetranslimagain}
\lim_{n\to\infty}\frac1n\log M_n(\lambda)=2\log(e^{\frac{\lambda}2}+1)-\log 4.
\end{equation}          
\end{proposition}

Let $I(x)$ denote the Legendre-Fenchel transform of the function appearing on the right hand side of \eqref{laplacetranslimagain}; that is,
\begin{equation}\label{Iagain}
I(x)=\sup_{\lambda\in\mathbb{R}}\Big(\lambda x-2\log(e^{\frac{\lambda}2}+1)+\log 4\Big),\ x\in\mathbb{R}.
\end{equation}
We will proof the following proposition.
\begin{proposition}\label{supL-Fagain}
The function $I$,  defined in \eqref{Iagain} and restricted to $x\in(0,1)$, is equal to $I^{\text{inc}}$ defined in \eqref{I-inc}.
\end{proposition}

\noindent \it Proof of Theorem \ref{incthm}.\rm\ In light of Propositions \ref{generatingpropagain} and \ref{supL-Fagain}, \eqref{ldinc} follows from the G\"artner-Ellis
theorem \cite{DZ}, except for the case $x=1$. The case $x=1$ is explained in Remark 2 following the statement of the theorem.
\hfill $\square$

We now turn to the proof of the two propositions.

\noindent \it Proof of Proposition \ref{generatingpropagain}.\rm\
In this paragraph we begin the proof for the case  $\eta=231$. In the next paragraph, we explain why the same proof works for $\eta=312$, and then continue to the end of the proof
for $\eta=231$. We have
$$
M_n(\lambda)=P_n^{\text{av}(231)}e^{\lambda L_n},\ \lambda\in\mathbb{R},\ n\ge0.
$$
 From the discussion in the first paragraph of the proof of Proposition \ref{generatingprop} and from the definition of $L_n$, it follows that
\begin{equation}\label{conddistagain}
\begin{aligned}
&L_n|\{\sigma_j=n\}\stackrel{\text{dist}}{=}L_{j-1}+L_{n-j},\ n\ge2,\ j=1,\cdots,n-1;\\
&L_n|\{\sigma_n=n\}\stackrel{\text{dist}}{=}L_{n-1}+1,\ n\ge2,
\end{aligned}
\end{equation}
where on the right hand side of \eqref{conddistagain}, for any $k$, the random variable $L_k$ is considered on $S_k$ under the measure $P_k^{\text{av}(231)}$, and
where on the right hand side of  the first line in \eqref{conddistagain},  $L_{j-1}$ and $L_{n-j}$ are independent.
From \eqref{conddistagain} and  \eqref{nprob}, we have
\begin{equation}\label{keyagain}
\begin{aligned}
&M_n(\lambda)=\sum_{j=1}^n P_n^{\text{av}(231)}(\sigma_j=n)E_n^{\text{av}(231)}(e^{\lambda L_n}|\sigma_j=n)=\\
&\sum_{j=1}^{n-1}\frac{C_{j-1}C_{n-j}}{C_n} E^{\text{av}(231)}_{j-1}e^{\lambda L_{j-1}}E^{\text{av}(231)}_{n-j}e^{\lambda L_{n-j}}+\frac{C_{n-1}}{C_n}
e^\lambda E^{\text{av}(231)}_{n-1}e^{\lambda L_{n-1}}=\\
&\sum_{j=1}^{n-1}\frac{C_{j-1}C_{n-j}}{C_n}M_{j-1}(\lambda)M_{n-j}(\lambda)+e^\lambda\frac{C_{n-1}}{C_n}M_{n-1}(\lambda),\ n\ge 2.
\end{aligned}
\end{equation}

In the case $\eta=312$, the same type of reasoning as in \eqref{nprob} shows
that $P_n^{\text{av}(312)}(\sigma_j=1)=\frac{C_{j-1}C_{n-j}}{C_n}$. Also, the same reasoning as in \eqref{conddistagain} gives
\begin{equation}\label{conddistagain312}
\begin{aligned}
&L_n|\{\sigma_j=1\}\stackrel{\text{dist}}{=}L_{j-1}+L_{n-j},\ n\ge2,\ j=2,\cdots,n;\\
&L_n|\{\sigma_1=1\}\stackrel{\text{dist}}{=}L_{n-1}+1,\ n\ge2,
\end{aligned}
\end{equation}
where on the right hand side of \eqref{conddistagain312}, for any $k$, the random variable $L_k$ is considered on $S_k$ under the measure $P_k^{\text{av}(312)}$, and
where on the right hand side of  the first line in \eqref{conddistagain312},  $L_{j-1}$ and $L_{n-j}$ are independent.
Using the these facts, one finds that the Laplace transform for this case also satisfies \eqref{keyagain}
Thus, it suffices to continue just for the case $\eta=231$.

Multiplying the leftmost and the rightmost expressions in \eqref{keyagain} by $C_nt^n$, we have 
\begin{equation}\label{fromkeyagain}
\begin{aligned}
&C_nM_n(\lambda)t^n= t\Big(\sum_{j=1}^{n-1}C_{j-1}M_{j-1}(\lambda)C_{n-j}M_{n-j}(\lambda)\Big)t^{n-1}+\\
&e^\lambda tC_{n-1}M_{n-1}(\lambda), \ n\ge2.
\end{aligned}
\end{equation}
Summing over $n$ gives
\begin{equation}\label{fromkeywithtagain}
\begin{aligned}
&\sum_{n=2}^\infty C_nM_n(\lambda)t^n= t\sum_{n=2}^\infty\Big(\sum_{j=1}^{n-1}C_{j-1}M_{j-1}(\lambda)C_{n-j}M_{n-j}(\lambda)\Big)t^{n-1}+\\
&e^\lambda t\sum_{n=2}^\infty C_{n-1}M_{n-1}(\lambda)t^{n-1}.
\end{aligned}
\end{equation}

Define
\begin{equation}\label{Gseriesagain}
G_\lambda(t)=\sum_{n=0}^\infty C_nM_n(\lambda)t^n.
\end{equation}
For use in some of the  calculations  below, note that $C_0=C_1=1$, $M_0(\lambda)=1$ and $M_1(\lambda)=e^\lambda$.
We have
\begin{equation}\label{G^2again}
G_\lambda^2(t)=\sum_{n=0}^\infty\Big(\sum_{k=0}^n C_kM_k(\lambda)C_{n-k}M_{n-k}(\lambda)\Big)t^n.
\end{equation}
The double sum on the right hand side of \eqref{fromkeywithtagain} can be written as
\begin{equation}\label{doublesumagain}
\begin{aligned}
&\sum_{n=2}^\infty\Big(\sum_{j=1}^{n-1}C_{j-1}M_{j-1}(\lambda)C_{n-j}M_{n-j}(\lambda)\Big)t^{n-1}=\\
&\sum_{n=2}^\infty\Big(\sum_{k=0}^{n-2}
C_kM_k(\lambda)C_{n-1-k}M_{n-1-k}(\lambda)\Big)t^{n-1}=\\
&\sum_{n=2}^\infty\Big(\sum_{k=0}^{n-1}
C_kM_k(\lambda)C_{n-1-k}M_{n-1-k}(\lambda)\Big)t^{n-1}-\sum_{n=2}^\infty C_{n-1}M_{n-1}(\lambda)t^{n-1}=\\
&\sum_{m=1}^\infty\Big(\sum_{k=0}^m
C_kM_k(\lambda)C_{m-k}M_{m-k}(\lambda)\Big)t^m-\sum_{m=1}^\infty C_mM_m(\lambda)t^m=\\
&(G^2_\lambda(t)-1)-(G_\lambda(t)-1)=G^2_\lambda(t)-G_\lambda(t),
\end{aligned}
\end{equation}
where \eqref{G^2again} has been used for the penultimate  equality.
From \eqref{fromkeywithtagain}, \eqref{doublesumagain} and the definition of $G_\lambda$, we obtain
$$
\big(G_\lambda(t)-e^\lambda t-1\big)
= t\big(G^2_\lambda(t)-G_\lambda(t)\big)+e^\lambda t\big(G_\lambda(t)-1\big),
$$
or equivalently
\begin{equation*}\label{Gequagain}
tG^2_\lambda(t)+\big((e^\lambda-1)t-1\big)G_\lambda(t)+1=0.
\end{equation*}
Since $G_\lambda(0)=1$, the quadratic formula yields
$$
G_\lambda(t)=\frac{1-(e^\lambda-1)t-\sqrt{\big((e^\lambda-1)t-1\big)^2-4t}}{2t},
$$
which we rewrite as
\begin{equation}\label{Gagain}
G_\lambda(t)=\frac{1-(e^\lambda-1)t-\sqrt{(e^\lambda-1)^2t^2-2(e^\lambda+1)t+1}}{2t}.
\end{equation}

From \eqref{Gagain}, it follows that the radius of convergence $R_\lambda$ of the power series representing $G_\lambda$ is the smaller of the absolute values of the two roots
of the quadratic polynomial 
$(e^\lambda-1)^2t^2-2(e^\lambda+1)t+1$. After a bit of algebra, one finds that the two roots are
$\frac{(e^{\frac\lambda2}\pm1)^2}{(e^\lambda-1)^2}$.
Thus, we have $R_\lambda=\frac{(e^{\frac\lambda2}-1)^2}{(e^\lambda-1)^2}$.
Consequently,
from \eqref{Gseriesagain}, we obtain
\begin{equation}\label{limsupagain}
\limsup_{n\to\infty}(C_nM_n(\lambda))^\frac1n=\frac{(e^\lambda-1)^2}{(e^{\frac\lambda2}-1)^2}.
\end{equation}
From \eqref{Catasymp} and \eqref{limsupagain}, we obtain
\begin{equation}\label{limsupMnagain}
\limsup_{n\to\infty}\frac1n\log M_n(\lambda)=2\log(e^\lambda-1)-2\log(e^{\frac\lambda2}-1)-\log 4=2\log(e^\frac\lambda2+1)-\log4.
\end{equation}
In light of \eqref{limsupMnagain}, to complete the proof of the proposition it suffices to show that the limit
$\lim_{n\to\infty}\frac1n\log M_n(\lambda)$ exists.
The proof of this is exactly the same as the corresponding proof in section \ref{pfaltthm}.
\hfill $\square$
\medskip

\noindent \it Proof of Proposition \ref{supL-Fagain}.\rm\ 
For $x\in(0,1)$, define
$$
g_x(\lambda)=\lambda x-2\log(e^\frac\lambda2+1)+\log 4, \ \lambda\in\mathbb{R}.
$$
We have
$$
g_x'(\lambda)=x-\frac{e^\frac\lambda2}{e^\frac\lambda2+1}.
$$
Setting the derivative equal to 0 and solving, we obtain
$\lambda=2\log\frac x{1-x}$.
Since $\lim_{\lambda\to\pm\infty}g_x(\lambda)=-\infty$, it follows that
$$
\begin{aligned}
&\sup_{\lambda\in\mathbb{R}}g_x(\lambda)=g_x(2\log\frac x{1-x})=2x\log x+2(1-x)\log(1-x)+\log4=\\
&2\big(x\log x+(1-x)\log(1-x)+\log 2\big)=I^{\text{inc}}(x),
\end{aligned}
$$
where $I^{\text{inc}}(x)$ is as in \eqref{I-inc}.
\hfill $\square$


\begin{thebibliography}{99}

\bibitem{AD}
Aldous, D. and  Diaconis, P.,
\emph{Longest increasing subsequences: from patience sorting to the Baik-Deift-Johansson theorem},
Bull. Amer. Math. Soc. (N.S.) \textbf{36} (1999),  413-432.

\bibitem{BDJ}
Baik, J.,  Deift, P. and  Johansson, K.,
\emph{On the distribution of the length of the longest increasing subsequence of random permutations},
J. Amer. Math. Soc. \textbf{12} (1999), 1119-1178.

\bibitem{B}  B\'ona, M., \emph{Combinatorics of Permutations},
Discrete Math. Appl. (Boca Raton)
CRC Press, Boca Raton, FL, (2012).


\bibitem{DZ} Dembo, A. and Zeitouni, O., \emph{Large Deviations Techniques and Applications},
Appl. Math. (N. Y.), 38,
Springer-Verlag, New York, (1998).


\bibitem{DZ99}
Deuschel, J.-D. and Zeitouni, O., 
\emph{On increasing sequences of IID samples},
Combin. Probab. Comput. \textbf{8} (1999), 247-263.



\bibitem{DHW}
Deutsch, E.,  Hildebrand, A. J. and Wilf, H. S.,
\emph{Longest increasing subsequences in pattern-restricted permutations}
Electron. J. Combin. \textbf{9} (2002),  Research paper 12, 8 pp.


\bibitem{FMW}
Firro, G., Mansour, T. and Wilson, M. C.,
\emph{Longest alternating subsequences in pattern-restricted permutations},
Electron. J. Combin. \textbf{14} (2007),  Research Paper 34, 17 pp.


\bibitem{LS}
Logan, B. F. and Shepp, L. A.,
\emph{A variational problem for random Young tableaux}
Advances in Math. \textbf{26} (1977), 206-222.

\bibitem{R}
Romik, D., \emph{The Surprising Mathematics of Longest Increasing Subsequences},
IMS Textb., \textbf{4}
Cambridge University Press, New York, (2015).

\bibitem{S} Sepp\"al\"ainen, T., \emph{Large deviations for increasing sequences in the plane},
Probab. Theory Relat. Fields \textbf{112} (1998), 221-244.  


\bibitem{SS} Simion, R. and Schmidt, F., \emph{Restricted permutations},
European J. Combin. \textbf{6} (1985),  383-406.


\bibitem{St} Stanely, R., \emph{Longest alternating subsequences of permutations}, Michigan Math. J. \textbf{57} (2008), 675-687.

\bibitem{VK} Veršik, A. M. and Kerov, S. V.,
\emph{Asymptotic behavior of the Plancherel measure of the symmetric group and the limit form of Young tableaux},
Dokl. Akad. Nauk SSSR \textbf{233} (1977), 1024-1027.

\bibitem{W}, Widom, H., \emph{On the limiting distribution for the longest alternating sequence
in a random permutation}, Electron. J. Combin. \textbf{13} (2006), Article R25.
\end{thebibliography}
\end{document}